
\input amstex
\documentstyle{amsppt}

\hcorrection{19mm}

\nologo
\NoBlackBoxes


\topmatter 

\title   Depth of pleated surfaces in toroidal cusps of hyperbolic
         3-manifolds
\endtitle
\author  Ying-Qing Wu$^1$
\endauthor
\leftheadtext{Ying-Qing Wu}
\rightheadtext{Depth of pleate surfaces in toroidal cusps}
\address Department of Mathematics, University of Iowa, Iowa City, IA
52242
\endaddress
\email  wu\@math.uiowa.edu
\endemail

\thanks  $^1$ Partially supported by NSF grant \#DMS 0203394
\endthanks

\abstract Let $F$ be a closed essential surface in a hyperbolic
3-manifold $M$ with a toroidal cusp $N$.  The depth of $F$ in $N$ is
the maximal distance from points of $F$ in $N$ to the boundary of $N$.
It will be shown that if $F$ is an essential pleated surface which is
not coannular to the boundary torus of $N$ then the depth of $F$ in
$N$ is bounded above by a constant depending only on the genus of
$F$.  The result is used to show that an immersed closed essential
surface in $M$ which is not coannular to the torus boundary components
of $M$ will remain essential in the Dehn filling manifold $M(\gamma)$
after excluding $C_g$ curves from each torus boundary component of
$M$, where $C_g$ is a constant depending only on the genus $g$ of the
surface.  \endabstract

\endtopmatter

\document

\define\proof{\demo{Proof}}
\define\endproof{\qed \enddemo}
\define\a{\alpha}
\redefine\b{\beta}

\define\r{\gamma}

\redefine\bdd{\partial}
\define\Int{\text{\rm Int}}
\define\BH{\Bbb H^3}
\input epsf.tex

\head  1.  Introduction  \endhead

In this paper a {\it surface\/} $F$ in a 3-manifold $M$ is a pair $F =
(S, \varphi)$, where $S$ is a connected, possibly nonorientable
surface, and $\varphi: S \to M$ is a continuous map which is an
immersion almost everywhere.  Concepts such as points on $F$ and
$\pi_1$-injectivity of $F$ are defined in the usual way.  See Section
2 for more details.

Let $N$ be a toroidal cusp in a hyperbolic 3-manifold $M$.  The {\it
depth\/} of a surface $F$ in $N$, denoted by $d_N(F)$, is defined as
the maximal distance from all points of $F \cap N$ to the boundary
torus $T$ of $N$.  When $F$ is disjoint from $N$, the depth of $F$ in
$N$ is defined as zero.  The depth of a closed surface $F$ in $N$ can
be arbitrarily large, as one can deform $F$ via homotopy to push $F$
deep into $N$.  On the other hand, if $F$ is closed, totally geodesic
and has a bounded genus, then the depth of $F$ in $N$ is bounded above
by a constant which is independent of $N$ and $M$.  See for example
[Ba1, Ba2].  Here it is important to assume that the genus of $F$ is
bounded, as the depth of closed totally geodesic surfaces in a cusp of
the figure 8 knot complement is unbounded.  (This follows from a
theorem of Leininger [Le] and the proof of Theorem 3.1 below.)

We would like to generalize this result to pleated surfaces, which
were invented by Thurston in [Th2].  Roughly speaking, a pleated
surface is a surface in a hyperbolic 3-manifold $M$ which is obtained
from a hyperbolic surface by bending along a geodesic lamination.
Because of that, it shares some properties of totally geodesic
surfaces; for example, it has the same area as the underlying
hyperbolic surface.  However, many important properties of totally
geodesic surfaces are not shared by pleated surfaces; for example, a
totally geodesic surface in $M$ is always $\pi_1$-injective, but that
is not true for pleated surfaces.

The following theorem gives an estimation of the depth of closed
essential pleated surfaces $F$ in a toroidal cusp $N$ of a hyperbolic
manifold $M$.  It shows that if $F$ is not coannular to the boundary
torus of a toroidal cusp $N$ then its depth in $N$ is bounded above by
a constant which depends only on the genus $g$ of $F$.

\proclaim{Theorem 1.1} Let $M$ be a complete hyperbolic 3-manifold,
and let $N$ be a toroidal cusp of $M$.  Let $F$ be an essential closed
pleated surface of genus $g$ in $M$ which is not coannular to the
torus $T = \bdd N$ in $M$.  Then $d_N(F) < \ln (2g\pi)$.  \endproclaim

The proof of this theorem will be given in Section 2.  One is referred
to that section for definitions of cusps and pleated surfaces.

We use the depth estimation of pleated surfaces to study the question
of how many Dehn fillings on a toroidal boundary component $T_0$ of a
compact hyperbolic manifold $M_0$ will preserve the
$\pi_1$-injectivity of an immersed essential surface $F$ in $M_0$.  A
surface $F = (S, \varphi)$ in $M_0$ is {\it incompressible\/} if it is
$\pi_1$-injective, i.e., $\varphi_*: \pi_1(S) \to \pi_1(M_0)$ is
injective.  $F$ is {\it essential\/} if it is incompressible and is
not homotopic to a surface on $\bdd M_0$.  $F$ is {\it coannular to
  $T_0$\/} if there is an incompressible annulus $A$ in $M_0$ with one
boundary component on each of $F$ and $T_0$, in which case the slope
of $A \cap T_0$ is called a coannular slope of $F$ on $T_0$.

The answer to this question is well understood if $F$ is embedded.
Denote by $M_0(\gamma)$ the manifold obtained from $M_0$ by Dehn
filling on $T_0$ along a slope $\gamma$.  It is
known that if $F$ is not coannular to $T_0$, then $F$ remains
essential in $M_0(\gamma)$ for all but at most three slopes $\gamma$
[Wu1], and if $F$ is coannular to $T_0$ with coannular slope
$\gamma_0$ then $F$ remains essential in $M_0(\gamma)$ unless
$\Delta(\gamma, \gamma_0) \leq 1$, where $\Delta (\alpha, \beta)$
denotes the geometric intersection number between the two slopes
[CGLS, Theorem 2.4.3].  Thus in the latter case $F$ remains essential
in $M_0(\gamma)$ for all but at most three lines of slopes $\gamma$.
Similar results hold for essential laminations and essential branched
surfaces.  See [Wu2].

For immersed essential surfaces, the above theorems can be generalized
in a weaker sense.  Suppose $F$ is an immersed essential surface in a
compact hyperbolic manifold $M_0$ with $\bdd M_0 = T_0$.  If $F$ is
not coannular to $T_0$, then there is a constant $K$ such that $F$
remains essential in $M_0(\gamma)$ for all but at most $K$ slopes
$\gamma$ [AR, Ba1].  If $F$ is coannular to $T_0$ with coannular
slopes $\alpha_1, ..., \alpha_n$, then there is a constant $K$ such
that $F$ remains essential in $M_0(\gamma)$ when $\Delta (\gamma,
\alpha_i) > K$ for all $i$ [Wu3, Theorem 5.3].

The major difference between the embedded and immersed cases is that
in the immersed case there is in general no upper bound for the
constant $K$ above.  It was shown in [Wu3] that for any constant $K$
there is an essential immersed surface $F$ in a hyperbolic manifold
$M_0$ with a coannular slope $\gamma_0$ on $T_0$ such that $F$ is
compressible in $M_0(\gamma)$ for all $\gamma$ satisfying
$\Delta(\gamma, \gamma_0) \leq K$.  For non-coannular case, Leininger
showed in [Le] that for any $K$ there is an essential surface in a
hyperbolic manifold $M_0$ which is totally geodesic (and hence is not
coannular to $T_0$), and is compressible in $M_0(\alpha)$ and
$M_0(\beta)$ with $\Delta (\alpha, \beta) > K$.  Thus there is no
universal upper bound on the geometric intersection number between
``killing slopes'' for $F$.  It is still an open question whether
there is a universal bound for the number of killings slopes for such
an $F$ in a hyperbolic manifold $M_0$, although we will show that no
such upper bound exists if $M_0$ is not assumed hyperbolic.  See
Theorem 3.7 below.

We use the result about depth of pleated surface in hyperbolic cusps
to study the question of whether there exists an upper bound for the
constant $K$ which depends only on the genus of the surface $F$.  The
following is a simplified version of Corrollary 3.2 in the special
case that the Dehn filling is performed on a single torus component
$T_0$ of $\bdd M_0$.

\proclaim{Theorem 1.2} Let $M_0$ be a compact orientable hyperbolic
3-manifold with $T_0$ a torus boundary component, and let $F$ be an
immersed essential surface in $M_0$ which is not coannular to $\bdd
M_0$.  Then there is a constant $C_1(g)$ ($<3200g^2$) depending only
on the genus $g$ of $F$, such that $F$ remains essential in
$M_0(\gamma)$ for all but at most $C_1(g)$ slopes $\gamma$ on $T_0$.
\endproclaim

The theorem shows that an upper bound for the constant $K$ above
depending only on the genus of $F$ does exist when $M_0$ is hyperbolic
and $F$ is not coannular to $T_0$.  This should be compared with
Theorem 3.5, which says that similar bound does not exist in the case
that $F$ is coannular to $T_0$.  It will also be shown (Theorem 3.7)
that the condition that $M_0$ be hyperbolic is crucial and cannot be
removed from the above theorem.  A general version of Theorem 1.2 (for
Dehn fillings on multiple boundary components of $M_0$) and its
variations are given in Theorem 3.1, Corollaries 3.2 and 3.3.

The constant $3200g$ in Theorem 1.2 is certainly not the best
possible.  A lower upper bound in the special case that $F$ is totally
geodesic can been found in [Ba1, Ba2].

\head  2.  Proof of Theorem 1.1
\endhead

Let $M_0$ be a compact, connected, orientable 3-manifold with $T$ a
set of tori on $\bdd M_0$, such that $M = M_0 - \bdd M_0$ admits a
complete hyperbolic structure.  Thus the universal covering space of
$M$ is the hyperbolic space $\Bbb H^3$.  We use the upper half space
model to identify $\BH$ with
$$\BH = \{(x,y,z) \,\, | \,\, z>0\}$$ which is endowed with a
riemannian metric $ds^2 = (1/z^2) dr^2$, where $dr^2$ is the Euclidean
metric of $\Bbb R^3$.

Let $\rho: \BH \to M$ be the covering map.  Recall that a {\it
horoball\/} of $\BH$ is a subset isometric to the set $\{(x,y,z) \,\,
| \,\, z \geq 1\}$.  The boundary of a horoball is called a {\it
horosphere.}  A {\it cusp\/} $N$ of $M$ at a torus boundary component
$T_0$ of $M_0$ is a regular neighborhood of $T_0$ in $M_0$ with $T_0$
removed, such that $\rho^{-1}(N)$ is a set of horoballs in $\BH$ with
disjoint interiors.  Each component of $\bdd M_0$ has a unique maximal
cusp in $M$, which has the property that there are a pair of horoballs
covering it which intersect each other at some point on their
boundary.  Given a cusp $N$, we may assume that $\rho: \BH \to M$ has
been chosen so that $\tilde N = \{(x,y,z)\,\, | \,\, z \geq 1 \}$ is a
horoball covering $N$.  Let $T = \bdd N$, and let $\tilde T = \bdd
\tilde N$.

Recall that a surface $F$ in $M$ is a pair $F = (S, \varphi)$, where
$S$ is a connected, possibly nonorientable surface, and $\varphi: S
\to M$ is a continuous map which is an immersion almost everywher.  A
surface $F$ in $M$ is incompressible if it is $\pi_1$-injective, i.e.,
$\varphi_*: \pi_1(S) \to \pi_1(M)$ is injective, in which case we
define the fundamental group of $F$ to be $\pi_1 F = \varphi_*
(\pi_1(S))$.  This is well defined up to conjugacy.  We can define
some other notions for $F$ in the obvious way.  For example a point on
$F$ is a pair $(x,\varphi)$ with $x\in S$, and a loop on $F$ is the
composition of $\varphi$ with a loop on $S$.  Thus if $F$ is
transverse to a surface $T$ in $M$ then $F \cap T$ is a set of loops
in $F$.  When there is no confusion, we will also refer to the set $F
= \varphi(S)$ as a surface.

A {\it rectifiable arc\/} in a Riemannian manifold $M$ is an arc
$\alpha$ which contains a nowhere dense closed subset $X$ such that
the closure of each component of $\alpha - X$ is a geodesic in $M$.
We use the definition given by Thurston in [Th2] for pleated surfaces:
A surface $F = (S, \varphi)$ in a hyperbolic manifold $M$ is a {\it
pleated surface\/} if $S$ is a complete hyperbolic surface, and
$\varphi: S \to M$ is a continuous map , such that 
\roster
\item[1] $\varphi$ is an isometry in the sense that every geodesic
segment in $S$ is taken to a rectifiable arc in $M$ which has the same
length, and

\item[2] for each point $x \in S$, there is at least one open
geodesic segment $l_x$ through $x$ which is mapped to a geodesic
segment in $M$.
\endroster

The {\it pleating locus\/} of $F$ is the set of points on $S$ which do
not have a neighborhood mapped by $\varphi$ to a totally geodesic disk
in $M$.  We refer the readers to [Th2] for some basic properties of
pleated surfaces.  In particular, it is known that the pleating locus
of a pleated surface $(S,\varphi)$ is a geodesic lamination on $S$,
which has measure $0$.

Now let $M$, $N$ and $F$ be as in Theorem 1.1.  We may assume that $F
= (S, \varphi)$ intersects $N$ nontrivially as otherwise we would have
$d_N(F) = 0$ and hence Theorem 1.1 holds.  Fix a point $p_0 =
\varphi(s_0) \in F \cap N$ such that the distance from $p_0$ to $T$
equals $d_N(F)$.  By definition and assumption $S$ is a closed
hyperbolic surface of genus $g$, so its universal covering $\tilde S$
is the hyperbolic plane $\Bbb H^2$.  Denote by $\rho_1$ the covering
map $\rho_1: \tilde S \to S$.  Let $\tilde s_0 \in \tilde S$ be a
point with $\rho_1(\tilde s_0) = s_0$, and let $\tilde \varphi: \tilde
S \to \tilde M$ be a lifting of $\varphi$ such that $\tilde p_0 =
\tilde \varphi(\tilde s_0)$ is a point in the horoball $\tilde N =
\{(x,y,z) \,\, | \,\, z \geq 1\}$, and the distance from $\tilde p_0$
to the horosphere $\tilde T = \bdd \tilde N$ equals $d_N(F)$.  We have
the following commutative diagram.
$$
\CD
(\tilde S, \tilde s_0)     @>\tilde \varphi   >>   (\BH, \tilde p_0) \\
@V{\rho_1}VV                          @VV{\rho}V \\
(S, s_0)     @>\varphi   >>   (M, p_0)
\endCD
$$
Note that $\tilde F = (\tilde S, \tilde \varphi)$ is a pleated surface
in $\tilde M$ with pleating locus $\rho_1^{-1}(\lambda)$, where $\lambda
\subset S$ is the pleating locus of $F$.

\proclaim{Lemma 2.1} Let $\tilde T_r = \{(x,y,z)\in \BH \,\,|\,\, z=r
\}$ be a horosphere in $\BH$ at level $r\geq 1$ which intersects
$\tilde F$ transversely.  Then each component of $\tilde F \cap \tilde
T_r$ is a circle bounding a disk $D$ in $\tilde F$ lying above $\tilde
T_r$ in $\BH$.  \endproclaim

\proof By assumption $F=(S, \varphi)$ intersects the torus $T_r =
\rho(\tilde T_r)$ transversely, so $\varphi^{-1}(T_r)$ is a set of
circles on $S$.  If some circle component $\r$ is mapped by $\varphi$
to an essential loop on $T_r$ then since $T_r$ is incompressible in
$M$, $\r$ must also be essential on $S$, so $F$ would be coannular to
$T_r$, contradicting the assumption of Theorem 1.1.  Therefore each
component of $\varphi^{-1} (T_r)$ is a circle which is trivial on $S$,
and hence bounds a disk $D'$ on $S$.  It follows that $\tilde \varphi
^{-1}(\tilde T_r)$ is a set of disjoint circles on $\tilde S$, each of
which bounds a disk $D$ in $\tilde S$.  We need to show that $\tilde
\varphi(D)$ lies above the horosphere $\tilde T_r$.

If this were not true, then since $D$ is compact, there is an {\it
interior\/} point $u$ of $D$ such that $\tilde \varphi(u)$ has minimal
$z$-coordinate in a neighborhood of $\tilde u$ in $D$.  On the other
hand, since $\tilde F = (\tilde S, \tilde \varphi)$ is a pleated
surface, there is a geodesic segment $\beta$ of $\BH$ lying on $\tilde
\varphi(\tilde S)$, containing the point $\tilde \varphi(u)$ in its
interior.  Since a geodesic in $\BH$ is either a euclidean half circle
perpendicular to the $xy$-plane or a vertical line parallel to the
$z$-axis, no point of $\beta$ can have locally minimal $z$-coordinate,
which is a contradiction.  \endproof

\proclaim{Lemma 2.2} Let $R(t)$ be the region on the $\Bbb H^2$ above
the horizontal line at height $1$ and inside the euclidean half circle
of radius $t>1$ centered at the origin, as shown in Figure 1.  Denote
by $A(R(t))$ the area of $R(t)$.  Then
$$A(R(t)) = 2 ( \sqrt{t^2 - 1} - \arctan \sqrt{t^2 - 1}) > 2t - 6.$$
\endproclaim

\proof  We have 
$$ \align 
A(R(t)) &= \iint_{R(t)} \frac 1{z^2} \, dy\, dz  
   = 2 \int_1^t \frac 1{z^2} \, dz \int_0^{\sqrt{t^2 - z^2}} dy \\
 & = 2 \int_1^t \frac{\sqrt{t^2 - z^2}}{z^2} \, dz \\
 & = 2\left[ \arctan \frac{\sqrt{t^2 - z^2}}z - \frac{\sqrt{t^2 -
   z^2}}z \right]_1^t \\ 
 & = 2 ( \sqrt{t^2 - 1} - \arctan \sqrt{t^2 - 1}) \\
 & > 2 ( (t-1) - \pi/2) > 2t - 6.  \qed
\endalign
$$
\enddemo

\bigskip
\leavevmode

\centerline{\epsfbox{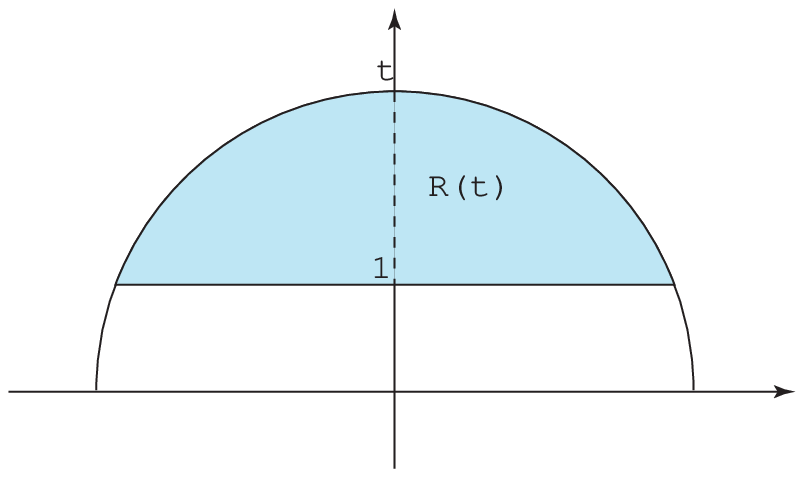}}
\nopagebreak
\bigskip
\centerline{Figure 1}
\bigskip

We identify $\Bbb H^2$ with the upper half $yz$-plane in $\Bbb H^3 =
\Bbb R^3_+$, and denote by $\psi: \BH \to \Bbb H^2$ the
euclidean orthogonal projection, i.e., 
$$ \psi(x,y,z) = (0, y,z).$$
Denote by $A(F)$ the area of a surface $F$ in $\BH$.  Clearly
$\psi$ is area decreasing, i.e., $A(\psi(F)) \leq A(F)$ for any
surface $F$ in $\Bbb H^3$.

\proclaim{Lemma 2.3}  Let $D$ be a disk in $\tilde N$ with $\bdd D
\subset \tilde T$.  Let $\alpha$ be a proper arc on $D$ such that
$\alpha' = \psi(\alpha)$ is a simple arc on $\Bbb H^2$.  Let
$\beta'$ be a horizontal arc on $\Bbb H^2$ joining the two endpoints of
$\alpha'$, and let $D'$ be the disk on $\Bbb H^2$ bounded by $\alpha' \cup
\beta'$.  Then $A(D) \geq 2 A(D')$.
\endproclaim

\proof Note that $\bdd D$ is mapped into the horizontal line at level
$z=1$ in $\Bbb H^2$, hence $\alpha' \cup \beta' = \partial D'$ is a
simple closed curve on $\Bbb H^2$.  The arc $\alpha$ cuts $D$ into
$D_1$ and $D_2$.  Note that $\psi$ maps $\bdd D_i$ to a curve on $\Bbb
H^2$ which is the union of $\alpha'$ and a path on the horizontal line
at $z=1$, hence if $u$ is a point in the interior of $D'$ then the
winding number of $\psi(\bdd D_i)$ around $u$ is $\pm 1$.  It follows
that $u \in \psi(D_i)$, in other words, $\psi(D_i)$ contains $D'$.
Since $\psi$ is area decreasing, $A(D') \leq A(\psi(D_i)) \leq
A(D_i)$, and the result follows.  \endproof

By considering the horosphere with $z$ coordinate $1+\epsilon$ if
necessary, we may assume that $\tilde F$ is transverse to the
horosphere $\tilde T$ at level $z=1$.  Let $D$ be the component of
$\tilde F \cap \tilde N$ containing the point $\tilde p_0 = (x_0, y_0,
z_0)$.  By definition of $\tilde p_0$, $z_0$ is maximal among
$z$-coordinates of all points of $\tilde F$.  By Lemma 2.1 $D$ is a
disk.  

\proclaim{Lemma 2.4} $A(D) \geq A(R(z_0)) > 2z_0 - 6.$ \endproclaim

\proof Without loss of generality we may assume that $\tilde p_0$ is
on the $z$-axis, so $\tilde p_0 = (0,0,z_0)$.  Recall that $\tilde F$
is a pleated surface.  If $\tilde p_0$ is on the pleating locus
$\tilde \lambda$ of $\tilde F$, then it is on a geodesic $\gamma$ of
$\BH$ contained in $\tilde F$.  If $\tilde p_0$ is not on the pleating
locus $\tilde \lambda$, then it lies in the interior of the closure of
a component of $\tilde F - \tilde \lambda$, which is a totally
geodesic surface $P$ in $\BH$ with boundary a disjoint union of
geodesics of $\BH$.  We will call $P$ an ideal polygon, although it
may have infinite area.

First assumee that $\tilde p_0$ is on a geodesic $\gamma$ of $\BH$
such that $\gamma \subset \tilde F$.  Rotating $\BH$ along the
$z$-axis if necessary, we may assume that $\gamma$ lies in the upper
half $yz$-plane, which is a hyperbolic plane $\Bbb H^2$.  Since
$\tilde p_0 =(0,0,z_0)$ has maximal $z$-coordinate on $\tilde F$ and
$\gamma \subset \tilde F$, the curve $\gamma$ is a euclidean half
circle of radius $z_0$ on $\Bbb H^2$ centered at the origin.  We now
apply Lemma 2.3 to the disk $D$ and the arc $\alpha = \gamma \cap D$.
Note that since $\gamma$ lies on the $yz$-plane, the projection
$\alpha' = \psi(\alpha)$ equals $\alpha$, which is the part of
$\gamma$ lying above the horizontal line $L$ at level $z=1$.  Thus the
disk $D'$ in Lemma 2.3 is exactly the region $R(z_0)$ defined in Lemma
2.2.  By Lemma 2.3 we have $A(D) \geq 2 A(D') = 2 A(R(z_0))$.

Now assume that $\tilde p_0$ is in an ideal polygon $P$ of $\BH$
contained in $\tilde F$.  Let $q$ be an ideal point of $P$, and let
$\gamma$ be a geodesic in $P$ passing through $\tilde p_0$ and having
one ideal endpoint at $q$.  Note that at least one component of
$\gamma - \tilde p_0$ lies in $P$, and $\gamma$ is either disjoint
from $\bdd P$ (and hence contained in $P$), or intersects $\bdd P$ at
a single point.  Up to a rotation of $\BH$ along the $z$-axis we may
assume that $\gamma$ lies on the $yz$-plane and the ideal point $q$ is
on the negative $y$-axis.  Let $R'$ be the region on $\Bbb H^2$ lying
on the left of the $z$-axis, below $\gamma$ and above the horizontal
line $z=1$.  Note that $A(R') = \frac 12 A(R(z_0))$.

If $\gamma \cap \bdd P = \emptyset$ then the proof follows as before.
So assume $\gamma$ intersects $\bdd P$ at a point $u$.  Then $u$ lies
on a geodesic $\beta$ on $\bdd P$, as shown in Figure 2.  Let
$\beta_1$ and $\beta_2$ be the two components of $\beta - u$, and let
$\alpha_i$ be the arc $(\gamma \cup \beta_i) \cap D$.  From Figure 2
it is clear that at least one of the projection of $\alpha_i$ on $\Bbb
H^2$ is an arc whose union with an arc on the horizontal line at $z=1$
bounds a disk $D'$ containing the region $R'$ above.  
By Lemma 2.3 it follows that $A(D) \geq 2 A(D') \geq 2 A(R') =
A(R(z_0)) > 2z_0 - 6$.  \endproof

\bigskip
\leavevmode

\centerline{\epsfbox{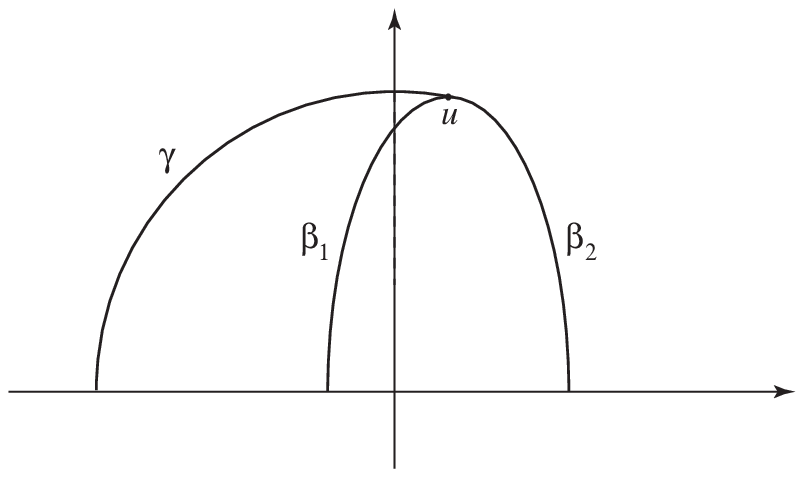}}
\nopagebreak
\bigskip
\centerline{Figure 2}
\bigskip

\demo{Proof of Theorem 1.1}  By Lemma 2.1 each component of
$\varphi^{-1}(N)$ is a disk on the hyperbolic surface $S$.  Since
$\rho_1: \tilde S \to S$ is a universal covering map, each component
of $\rho_1^{-1}(\varphi^{-1}(N))$ is a disk $D$ which maps
homeomorphically by $\rho_1$ to a component $D'$ of $\varphi^{-1}(N)$;
therefore
$$A(D) = A(D') \leq A(S) = 2\pi |\chi(S)| = 2 \pi (2g - 2).$$
Since $\tilde \varphi: \tilde S \to \BH$ is a pleated surface, the
area of $\tilde \varphi(D)$ equals the area of $D$.
Thus by Lemma 2.4 we have
$$2 \pi (2g-2) \geq A(D') = A(D) > 2 z_0 - 6,$$ where $z_0$ is the
maximum $z$-coordinate of points on $\tilde \varphi(\tilde S)$.  Thus
$z_0 < 2g\pi$.  Now the hyperbolic distance from a point $p=(x,y,z)
\in \BH$ to the horosphere $\tilde T$ at level $z=1$ is $|\ln z|$.
Therefore the depth of $F$ in $N$ satisfies $ d_N(F) = \ln z_0 < \ln
(2g\pi).$ \endproof

\head 3.  Essential surfaces and Dehn filling
\endhead

Let $M_0$ be a compact orientable hyperbollic 3-manifold with $\bdd
M_0$ a set of tori, and let $M = M_0 - \bdd M_0$.  A {\it maximal set
  of cusps\/} in $M$ is a set of cusps $N$ in $M$, one for each end of
$M$, such that (i) the interiors of the cusps are disjoint, and (ii)
$N$ is maximal subject to (i).  Note that $N$ may not be unique if
$M_0$ has more than one boundary components.

Let $N$ be a set of cusps in $M$, one for each end of $M$.  Let $T =
\bdd N$, which is a set of tori in $M$, and is parallel to $\bdd M_0$
in $M_0$.  The hyperbolic structure of $M$ induces a Euclidean metric
on $T$.  The {\it $T$-length\/} $l_T(\gamma)$ of a curve $\gamma$ on
$\bdd M_0$ is defined to be the length of a geodesic on $T$ homotopic
to $\gamma$ in $M_0$.  It is known that if $N$ is maximal and $T =
\bdd N$ then $l_T(\gamma) \geq 1$ for any nontrivial curve $\gamma$ on
$\bdd N$.  See [BH].

A {\it multiple slope\/} is a set of essential simple closed curves
$\gamma = (\gamma_1, ..., \gamma_k)$, such that there is at most one
$\gamma_i$ for each boundary component of $M_0$.  Define $M_0(\gamma)$
to be the manifold obtained by attaching $k$ solid tori to $M_0$ such
that $\gamma_i$ bounds a meridian disk in the $i$-th attached solid
torus.  

\proclaim{Theorem 3.1} Let $F = (S,\varphi)$ be a closed
incompressible surface of genus $g$ in a compact, orientable,
hyperbolic 3-manifold $M_0$ with $\partial M_0$ a set of tori, and
assume that $F$ is not coannular to $\bdd M_0$.  Then there is a
constant $C_0(g) < 4\pi^2 g$ depending only on $g$, such that if $N$
is a maximal set of cusps of $M = M_0 - \partial M_0$, and $\gamma =
\gamma_1 \cup ... \cup \gamma_k$ is a multiple slope on $\bdd M_0$
such that $l_T(\gamma_i) > C_0(g)$ for all $i$, then $F$ remains
incompressible in $M_0(\gamma)$.  \endproclaim

\proof This follows from Theorem 1.1 and a theorem of Bart [Ba1,
Theorem 1.1], which says $F$ remains essential in $M(\gamma)$ if one
can choose a set of cusps disjoint from the projection of the convex
hull of the limit set of $\pi_1(F)$, and the length of each $\gamma_i$
is greater than $2\pi$.  We give a proof for the convenience of the
readers.

By assumption $M = M_0 - \bdd M_0$ admits a complete hyperbolic
structure of finite volume.  Recall that the group of the surface $F =
(S,\varphi)$ is defined to be the subgroup $G = \varphi_*(\pi_1(S))$
of $\pi_1(M)$.  By Thurston [Th1], either $F$ is a virtual fiber, or
$G$ is geometrically finite; but since $F$ is closed while $M$ is
noncompact, $F$ cannot be a virtual fiber.  Hence it must be
geometrically finite.

Let $H(G)$ be the convex hull of the limit set of $G$ in $\BH$.
Consider the covering space $M_1$ of $M$ corresponding to the subgroup
$G$.  Let $\rho': \BH \to M_1$ and $\rho'': M_1 \to M$ be the covering
maps.  Then $\rho'(H(G))$ is a compact set in $M_1$ homeomorphic to $S
\times I$ (or $S$ if $G$ is Fuchsian), and hence $\rho(H(G)) = \rho''
(\rho'(H(G)))$ is compact in $M$.  See [Mg].  The boundary of
$\rho'(H(G))$ consists of two pleated surfaces $F'_1, F'_2$ in $M_1$
($F'_1 = F'_2$ if $G$ is Fuchsian), which map to pleated surfaces
$F_1, F_2$ in $M$ by $\rho''$.  Since $\rho: \BH \to M$ is a covering
map, its restrictioin on $H(G)$ is an immersion; hence the frontier of
$\rho(H(G))$ is contained in $\rho(\bdd H(G)) = F_1 \cup F_2$.  By
Theorem 1.1, the depth of $F_i$ in a cusp of $M$ is bounded by $\ln
(2g\pi)$, so the depth of $\rho(H(G))$ in $N$ is also bounded by
$\ln(2g\pi)$.

For simplicity of notations let us assume $k=1$.  The proof for the
general case is similar.  Thus $N$ is a cusp adjacent to a torus
component $T_0$ of $\bdd M_0$ containing the Dehn filling slope
$\gamma$, and the length of $\gamma$ on $T = \bdd N$ is greater than
$4g\pi^2$.  Let $T'$ be the euclidean torus in the cusp $N$ at depth
$\ln(2g\pi)$.  By lifting to the universal covering space $\BH$ one
can see that the length of a geodesic curve $\gamma'$ on $T'$ isotopic
to $\gamma$ in $N$ is greater than $(4g\pi^2)/(2g\pi) = 2\pi$.

Denote by $M'$ the manifold obtained from $M$ by cutting off the cusp
bounded by $T'$.  By the $2\pi$ theorem of Gromov-Thurston [GT, BH],
if $\gamma$ is a slope on $T'$ with geodesic length greater then
$2\pi$ then the hyperbolic structure on $M'$ can be extended over the
Dehn filling solid torus to obtain a negatively curved metric on the
Dehn filled manifold $M'(\gamma)$, which is clearly homeomorphic to
$M_0(\gamma)$ by a homeomorphism rel $\rho(H(G))$.  Up to homotopy we
may assume that $F$ is in $\rho(H(G))$.  If $\alpha$ is an essential
closed curve on $S$, then $\varphi(\alpha)$ is homotopic to a geodesic
$\a'$ in $M$, which lifts to a geodesic in the convex hull $H(G)$, and
hence $\a'$ lies in $\rho(H(G))$, which is a subset of $M'$.  It
follows that $\a'$ remains a geodesic in $M'(\r)$, which must be
essential because $M'(\r)$ is negatively curved.  Therefore $F =
(S,\varphi)$ is $\pi_1$-injective in $M'(\gamma)$, and hence in
$M_0(\gamma)$.  \endproof

Given a euclidean torus $T$, there are only finitely many curves of
length at most $t$.  In fact, there is a constant $C(t)$, independent
of $T$, such that if the minimum length of closed geodesics on $T$ is
$1$ then there are at most $C(t)$ simple closed geodesics $\r$ on $T$
with $l_T(\r) \leq t$.  Using the idea in the proof of [BH, Lemma 12]
one can show that
$$ C(t) \leq \frac {\frac {\pi}{2} (t+\frac 12)^2} {\frac {\sqrt
    3}{2}} + 1, $$
so for $t\geq 4\pi^2$ we have $C(t) < 2t^2$.  The following result is
now an immediate corollary of Theorem 3.1, as it is known that the
minimal length of a closed geodesic on the boundary of a maximal set
of cusps is at least 1 [BH].

\proclaim{Corollary 3.2} Let $M_0$ be a compact orientable hyperbolic
3-manifold with toroidal boundary, and let $F$ be an immersed
essential surface of genus $g$ in $M_0$ which is not coannular to
$\bdd M_0$.  Then after excluding at most $C(4\pi^2 g)$ ($< 3200g^2$)
curves on each component of $\bdd M_0$, $F$ remains essential in all
$M(\gamma)$.  \endproclaim

Denote by $\Delta (\gamma_1, \gamma_2)$ the minimal geometric
intersection number between two slopes $\gamma_1$ and $\gamma_2$ on a
torus $T$.  There is another constant $C'(t)$ such that if (i) $T$ is
a euclidean torus such that the minimal length of closed geodesics on
$T$ is at least 1, and (ii) $\gamma_1$ and $\gamma_2$ are closed
geodesics of length at most $t$ on $T$, then $\Delta(\gamma_1,
\gamma_2) < C'(t)$.  Using the idea in the proof of [BH, Theorem 14]
one can show that 
$$ C'(t) < \frac 2{\sqrt 3} t^2 < 2t^2.$$

\proclaim{Corollary 3.3} Let $M_0$ be a compact orientable hyperbolic
3-manifold with toroidal boundary, and let $F$ be an immersed
essential surface of genus $g$ in $M_0$ which is not coannular to
$\bdd M_0$.  Let $\gamma_1, \gamma_2$ be slopes on a torus component
$T_0$ of $\bdd M_0$.  If $F$ is compressible in both $M_0(\gamma_1)$
and $M_0(\gamma_2)$, then $\Delta (\gamma_1, \gamma_2) < C'(4 \pi^2 g)
< 3200g^2$.  \endproclaim

This should be compared with Leininger's result [Le], which says that
there is no universal bound for $\Delta(\gamma_1, \gamma_2)$ if
there is no constraint about the genus of $F$.

In Theorem 3.1 we assumed that the manifold $M_0$ is hyperbolic, and
the surface $F$ is not coannular to $\bdd M_0$.  Neither condition can
be removed.  The following theorem is in [Wu3, Theorem 1.1].

\proclaim{Theorem 3.4} Let $M$ be a compact, orientable, hyperbolic
3-manifold, and let $T$ be a torus component of $\bdd M$.  Let $F$ be
a closed essential surface in $M$.  

(1) If $F$ is coannular to $T$ with coannular slopes $\r_1, ...,
\r_n$, then there is an integer $K$ such that $F$ remains
$\pi_1$-injective in $M(\r)$ for all $\r$ satisfying $\Delta(\r, \r_i)
\geq K$, $i=1,...,n$.

(2) There is no universal bound on the constant $K$ in (1): For any
constant $C$, there is a closed essential surface $F$ in a hyperbolic
manifold $M$ with coannular slope $\b$ on $T=\bdd M$, such that $F$ is
not $\pi_1$-injective in $M(\r)$ for all $\r$ with $\Delta(\r, \b)
\leq C$.  
\endproclaim

Part (1) of the theorem is a finiteness theorem, while part (2) shows
there is no upper bound in general for the number of lines of
surgeries which will kill the essential surface.  The genera of the
surfaces used in the proof of part (2) increase as $C$ increases.
Because of Theorem 3.1, one may wonder if a similar theorem would be
true for surfaces coannular to $\bdd M$.  More precisely, does there
exist a constant $C(g)$ depending only on the genus of $F$, such that
$F$ remains essential in $M(\gamma)$ for all but at most $C(g)$ lines
of slopes $\gamma$?  The following theorem shows that no such constant
$C(g)$ exist.  Recall that with respect to the standard
meridian-longitude pair $(m,l)$ of a knot $K$ in $S^3$, a slope
$\gamma = pm+ql$ on $\bdd N(K)$ is represented by a rational number
$p/q$ or $1/0$.  See [Ro].  Thus $M(\gamma) = M(p/q)$.

\proclaim{Theorem 3.5} Let $K$ be a non-fibred hyperbolic knot in
$S^3$ with genus $g$, and let $M = S^3 - \Int N(K)$ be the knot
exterior.  Then for any constant $C$ there is an immersed essential
surface $F$ of genus $2g$ with the longitude of $K$ as a coannular
slope, such that $F$ is compressible in $M(\gamma)$ for all $\gamma =
p/q$ satisfying $0<p\leq C$.  \endproclaim

\proof Let $S_1$ be a minimal Seifert surface of $K$ in $M$.  Let $S$
be the double of $S_1$, which is an abstract surface obtained by
taking two copies of $S_1$ joined by an annulus $A$.  Let $\varphi_k:
S \to M$ be the map which sends the two copies of $S_1$ to $S_1$, and
$A$ to an annulus wrapping $k$ times around the torus $T = \bdd M$.
The surface $F_k = (S, \varphi_k)$ is called a Freedman tubing of
$S_1$ (with wrapping number $k$).  See [FF, CL, Wu3] for more details.
Since $K$ is hyperbolic and not fibred, by a theorem of
Freedman-Cooper-Long [FF, CL, Li], $F_k$ is incompressible in $M$ when
$k$ is sufficiently large.  See [Wu3] for an alternative proof.

\proclaim{Lemma 3.6} If $k$ is a multiple of $p$ then the surface $F_k$
is compressible in $M(p/q)$.  \endproclaim

\proof Put $k = pr$.  Let $\a$ be an essential arc on the annulus
$\varphi_k(A)$ which wraps $k$ times around the meridian $m$ and $0$
times along the longitude $l$.  Then $\a \cdot (l^{qr})$ represents
$km + (qr)l = r(pm + ql)$ in $\pi_1(T)$, and hence is null-homotopic
in the Dehn filling solid torus $V$.  Thus $\a$ is homotopic to
$l^{-qr}$ in $V$, so the annulus $\varphi_k(A)$ is rel $\partial A$
homotopic to a degenerate annulus with image on the boundary of the
Seifert surface $S_1$.  Therefore the immersed surface $F_k$ is
homotopic in $M(p/q)$ to an embedding of $S$ in $M$ with image $\bdd
N(S_1)$, the boundary of a regular neighborhood of $S_1$.  Since
$N(S_1)$ is a handlebody, $\bdd N(S_1)$ is compressible in $M$, hence
$\varphi_k$ is also compressible in $M(p/q)$.  \endproof

Now given any positive integer $C$, let $k$ be a positive integer
which is a multiple of all $p$ with $0<p\leq C$.  We may assume that
$C$ is sufficiently large, so the surface $F_k$ is essential by the
Freedman-Cooper-Long Theorem.  By Lemma 3.6, the surface $F_k$ is
compressible in $M(p/q)$ for all $p/q$ with $0<p\leq C$.  This
completes the proof of Theorem 3.5.  \endproof

The following theorem shows that Theorem 3.1 is false if the manifold
$M$ is not assumed to be hyperbolic.

\proclaim{Theorem 3.7} There exists a non-hyperbolic compact
orientable irreducible 3-manifold $M$ with $\bdd M = T$ a torus, a
fixed closed connected orientable surface $S$ of positive genus, and
an immersion $\varphi_{n(k)}: S \to M$ for any positive integer $k$,
such that the essential surfaces $F_k = (S, \varphi_{n(k)})$ satisfy
the following conditions.

(i)  $F_k$ is incompressible;

(ii) $F_k$ is not coannular to $T$;

(iii) There are $k$ distinct slopes $\r_1, ..., \r_k$ on $T$ such that
$F_k$ is compressible in $M(\r_i)$ for $i=1, ..., k$.
\endproclaim

\proof Let $K_1$ be a hyperbolic knot in $S^3$ which is not a fibred
knot.  Let $K = C_{p,q}(K_1)$ be a $(p,q)$ cable knot of $K_1$, and let
$M = S^3 - \Int N(K)$ be the exterior of $K$.  Denote by $V = N(K_1)$ a
regular neighborhood of $K_1$ which contains $K$ as a cable knot.  Let
$M_1 = S^3 - \Int V$, and let $X = V - \Int N(K)$.  Let $T_1 = \bdd
M_1$.  Thus $M = M_1  \cup_{T_1} X$.  

To construct the immersed surfaces $F_k$, let $S_1$ be a minimal genus
Seifert surface of the hyperbolic knot $K_1$, and let $\varphi_{n}:
S \to M_1$ be a Freedman tubing of two copies of $S_1$ with wrapping
number $n$.  By the proof of Theorem 3.5 we may choose $n(k)$ so
that (a) $F_k = (S, \varphi_{n(k)} )$ is incompressible, and (b) $F_k$
is compressible in $M_1(\gamma)$ for all $\gamma=r/s$ satisfying
$1\leq r \leq kpq+1$.  In particular, $F_k$ satisfies condition (i).

Let $\r_i = (ipq+1)/i$.  By a theorem of Gordon [Go, Corollary 7.3],
we have
$$M(\r_i) = M(\frac{ipq + 1}{i}) = M_1(\frac{ipq+1}{iq^2}).$$ 
By definition, $F_k$ is compressible in $M_1((ipq+1)/iq^2)$ for all
$i$ satisfying $1\leq i \leq k$.  Therefore $F_k$ satisfies condition
(iii).  It remains to show that $F_k$ is not coannular to $T$.  

First notice that the only coannular slope of $F_k$ in the manifold
$M_1$ is the longitude.  This is because $F_k$ lifts to a closed
surface in the universal abelian covering of $M_1$, so any closed
curve on $F_k$ lifts, but the longitude is the only primitive closed
curve on $T_1 = \bdd M_1$ that lifts.

Now suppose $A$ is an essential annulus in $M$ with one boundary
component in each of $F$ and $T$.  By a homotopy we may assume that
$A$ is transverse to $T_1$, and the number of components of $A \cap
T_1$ is minimal.  Since $T_1$ is incompressible in $M$, each component
of $A \cap M_1$ and $A \cap X$ is an essential annulus in $M_1$ or
$X$, except the one which has a boundary component on $F$, denoted by
$A_1$.  By the above, $\beta = A_1 \cap T_1$ is a longitude of $K_1$.
On the other hand, the component $A_2$ of $A \cap X$ which has $\beta$
as a boundary component is an essential annulus in $X$, and it is well
known that such an annulus is isotopic to a union of fibers [Ja,
Theorem VI.34]; in particular, the longitude $\beta$ of $K_1$ is a
fiber of $X$.  This is a contradiction because after trivial filling
along the meridional slope $m$ on $T$ the Seifert fibration on $X$
extends to a Seifert fibration of the solid torus $X(m)$ containing
$K$ as a regular fiber, so a fiber on $T_1$ has slope $p/q$.
\endproof

\enddocument